\documentclass[a4paper,10pt]{amsart}
\usepackage[english]{babel}
\input xy
\xyoption{all}
\usepackage{amssymb,amsmath,amsthm}

\newcommand{\R}{\mathbb{R}}

\newcommand{\lie}[1]{\mathfrak{#1}}     
\newcommand{\Z}{\mathbb{Z}}

\newcommand{\hook}{\lrcorner\,}
\newcommand{\LieG}[1]{\mathrm{#1}}      

\newcommand{\Spin}{\mathrm{Spin}}
\newcommand{\SU}{\mathrm{SU}}
\newcommand{\PSU}{\mathrm{PSU}}
\newcommand{\SO}{\mathrm{SO}}

\newcommand{\Gtwo}{{\mathrm{G}_2}}
\newcommand{\so}{\mathfrak{so}}
\newcommand{\gl}{\mathfrak{gl}}
\newcommand{\su}{\mathfrak{su}}

\newcommand{\GL}{\mathrm{GL}}

\newcommand{\dfn}[1]{\emph{#1}}
\newcommand{\leftdoublebracket}{[\![}
\newcommand{\rightdoublebracket}{]\!]}

\DeclareMathOperator{\Hom}{Hom}

\DeclareMathOperator{\codim}{codim}

\theoremstyle{plain}
\newtheorem{proposition}{Proposition}
\newtheorem{theorem}[proposition]{Theorem}
\newtheorem{lemma}[proposition]{Lemma}
\newtheorem{corollary}[proposition]{Corollary}
\theoremstyle{definition}
\newtheorem{definition}[proposition]{Definition}

\theoremstyle{remark}
\newtheorem*{remark}{Remark}
\newcommand{\Span}[1]{\operatorname{Span}\left\{#1\right\}}

\begin{document}
\title{Embedding into manifolds with torsion}
\author{Diego Conti}
\address{Dipartimento di Matematica e Applicazioni, Universit\`a di Milano Bicocca, via Cozzi 53, 20125 Milano, Italy.}
\email{diego.conti@unimib.it}
\thanks{This research has been carried out despite the effects of the Italian law 133/08. This law drastically reduces public funds to public Italian universities, with dangerous consequences for free scientific research, and will prevent young researchers from obtaining a position, either temporary or tenured, in Italy. The author is protesting against this law to obtain its cancellation (see {\tt http://groups.google.it/group/scienceaction}).}

\maketitle
\sloppy
\begin{abstract}
We introduce a class of special geometries associated to the choice of a differential graded algebra contained in $\Lambda^*\R^n$. We generalize some known embedding results, that effectively characterize the real analytic Riemannian manifolds that can be realized as submanifolds of a Riemannian manifold with special holonomy, to this more general context. In particular, we consider the case of hypersurfaces inside nearly-K\"ahler and $\alpha$-Einstein-Sasaki manifolds, proving that the corresponding evolution equations always admit a solution in the real analytic case.
\end{abstract}
There are a number of known results concerning the problem of embedding a generic Riemannian manifold into a manifold with special holonomy. At the very least for this problem to make sense the embedding should be isometric, but often extra conditions are required. A classical example is that of special Lagrangian submanifolds of Calabi-Yau manifolds, introduced in \cite{HarveyLawson:CalibratedGeometries} in the context of minimal submanifolds; it is a result of \cite{Bryant:Calibrated} that every compact, oriented real analytic Riemannian $3$-manifold can be embedded isometrically as a special Lagrangian submanifold in a $6$-manifold with holonomy contained in $\SU(3)$ (see also \cite{Matessi:Isometric} for a generalization). A similar result holds for coassociative submanifolds in $7$-manifolds with holonomy contained in $\Gtwo$ (see \cite{Bryant:Calibrated}); the case of codimension one embeddings in manifolds with special holonomy was considered in \cite{Hitchin:StableForms, ContiSalamon}.

Whilst the above results have been obtained using characterizations in terms of differential forms, the codimension one embedding problem can be uniformly rephrased using the language of spinors. Indeed, if $M$ is a Riemannian spin manifold with a parallel spinor, then a hypersurface inherits a generalized Killing spinor $\psi$, namely satisfying
\[\nabla_X\psi = \frac12 A(X)\cdot\psi\]
where $A$ is a section of the bundle of symmetric endomorphism of $TM$, corresponding to the Weingarten tensor, and the dot represents Clifford multiplication. One can ask if, given a Riemannian spin manifold $(N,g)$ with a generalized Killing spinor $\psi$, the spinor can be extended to a parallel spinor on a Riemannian manifold $M\supset N$ containing $N$ as a hypersurface; if so, $M$ is Ricci flat and the second fundamental form is determined by the tensor $A$, which is an intrinsic property of $(N,g,\psi)$. Counterexamples exist in the smooth category \cite{Bryant:NonEmbedding}; in the real analytic category, the general existence of such an embedding has been established in \cite{BarGauduchonMoroianu} under the assumption that the tensor $A$ is Codazzi, and in \cite{ContiSalamon} for six-dimensional $N$. A characterization of the $(N,g,\psi)$ as above that embed into an irreducible $M$ follows from Theorem~\ref{thm:AbstractEmbedding} in this paper. The existence of an embedding for any real analytic $(N,g,\psi)$ has been recently proved by Ammann \cite{Ammann}.

\smallskip
The spinor formulation lends itself to a natural variation, where one replaces the parallel spinor on the ambient manifold with a real Killing spinor, namely a spinor $\psi$ for which $\nabla_X\psi = \frac12 X\cdot\psi$. Compact manifolds with a Killing spinor are classified in \cite{Bar}, and they are either Einstein-Sasaki (possibly $3$-Einstein-Sasaki) manifolds, nearly-K\"ahler six-manifolds, nearly-parallel $\Gtwo$-manifolds, or round spheres. Aside from the somewhat degenerate case of the sphere, each of these geometries can be associated to a finite-dimensional differential graded algebra (DGA). For instance, nearly-K\"ahler geometry corresponds to the algebra of $\SU(3)$-invariant elements in the exterior algebra $\Lambda^*\R^6$, which is generated by a $2$-form $\omega$ and two $3$-forms $\psi^\pm$; the differential graded algebra structure is given by the operator
\[d(\omega)=3\psi^+, \quad d(\psi^-)=-2\omega\wedge\omega, \quad d(\psi^+)=0.\]
A hypersurface inside a manifold with a nearly-K\"ahler $6$-manifold has an induced structure, called a nearly hypo $\SU(2)$-structure, and the codimension one embedding problem for nearly hypo $5$-manifolds amounts to determining whether the ``nearly hypo evolution equations'' admit a solution \cite{Fernandez:NearlyHypo}.

In this paper we give a characterization of the DGA's that define a special geometry (Lemma~\ref{lemma:ZEmpty}), and classify these DGA structures on the algebras $(\Lambda^*\R^{2n})^{\SU(n)}$ and $(\Lambda^*\R^{2n+1})^{\SU(n)}$ (Propositions~\ref{prop:derivation} and \ref{prop:oddderivation}). We define the exterior differential system associated to such a DGA, and provide examples in which it is involutive (Theorem~\ref{thm:involutive}); as a byproduct, we obtain that the intrinsic torsion of an $\SU(n)$\nobreakdash-structure on a manifold of dimension $2n+1$ is entirely determined by the exterior derivative of the defining forms, like in the well-known even-dimensional case (Proposition~\ref{prop:SUnIsStronglyAdmissible}). We formulate the embedding problem for geometries defined by arbitrary DGA's and submanifolds of arbitrary codimension, and determine a sufficient condition for the embedding to exist by means of Cartan-K\"ahler theory (Theorem~\ref{thm:AbstractEmbedding}). This result requires that the manifold and  metric be real analytic, which appears to be a necessary restriction due to the counterexamples of \cite{Bryant:NonEmbedding}.
Though the sufficient condition is not very explicit, we show it to hold in some explicit examples, including those corresponding to an ambient geometry with a real Killing spinor (Theorems~\ref{thm:nearlyhypo} and \ref{thm:go}), and those where the DGA is generated by a single stable form (Proposition~\ref{prop:stable}). In addition, we study a new geometry defined by a $5$-form in nine dimensions whose stabilizer $\SO(3)\subset\SO(9)$ acts irreducibly on $\R^9$, proving an embedding result in that case (Theorem~\ref{thm:so3Embedding}). A different irreducible representation of $\SO(3)$ has been studied in \cite{Nurowski:SO3}, along with its associated five-dimensional geometry, which is however defined by a tensor rather than differential forms. In seven dimensions, the $\SO(3)$-invariant forms are also invariant under $\Gtwo$; the nine-dimensional case seems more in line with the DGA approach of this paper, since the irreducible action of $\SO(3)$ on $\R^9$ preserves two differential forms, each with stabilizer $\SO(3)$.

\section{$\PSU(3)$-structures and calibrated submanifolds}
\label{sec:psu3}
In this section we give a concrete introduction to the embedding problem, in the context of $\PSU(3)$-structures on $8$-manifolds. This type of structure has been studied in \cite{Hitchin:StableForms, Witt:SpecialMetrics}.

The structure group $\PSU(3)=\SU(3)/\Z_3$ acts faithfully on $\R^8$ via the adjoint representation; a $\PSU(3)$-structure on an $8$-manifold $M$ is a reduction of the bundle of frames to the group $\PSU(3)\subset\SO(8)$. The group $\PSU(3)$ fixes a $3$-form on $\R^8$, which in terms of an orthonormal basis $e^1,\dotsc, e^8$ can be written as
\[\rho=e^{123}+\frac12 e^1\wedge (e^{47}-e^{56})+\frac12 e^1\wedge (e^{46}+e^{57})+\frac12 e^3\wedge (e^{45}-e^{67})+\frac{\sqrt3}2e^8\wedge (e^{45}+e^{67})\;.\]
This form is stable in the sense of \cite{Hitchin:StableForms}, a fact which plays a significant r\^ole in the proof of Theorem~\ref{thm:psu3Embedding} below.

By construction, an $8$-manifold $M$ with a $\PSU(3)$-structure carries two canonical forms, which we also denote by  $\rho$ and $*\rho$. Requiring these forms to be parallel is too strong a condition in order to obtain interesting examples (see \cite{Hitchin:StableForms}), and one is led to consider  the weaker conditions
\begin{equation}
 \label{eqn:PSU3Geometries} d\rho=0, \quad d*\rho=0, \quad d\rho=0=d*\rho.
\end{equation}
We shall see in Section~\ref{sec:Preparatory} that to each geometry in \eqref{eqn:PSU3Geometries} one can associate an exterior differential system, which is involutive in the first two cases (see Proposition~\ref{prop:stable}), though it does not appear to be involutive in the last case.

Up to a normalization constant, we can assume that the comass of $\rho$, namely the maximum of the map
\begin{equation}
\label{eqn:comass}
Gr_3^+(\R^8)\to\R, \quad E\to \rho(E),
\end{equation}
is equal to one. Then if $M$ is an $8$-manifold with a $\PSU(3)$-structure such that the associated form $\rho$ is closed, $\rho$ is a calibration in the sense of \cite{HarveyLawson:CalibratedGeometries}. An oriented three-dimensional $N\subset M$ is said to be calibrated by $\rho$ if $\rho(T_xN)=1$ for all $x$ in $N$. In analogy with \cite{Bryant:Calibrated}, we can prove the following:
\begin{theorem}
\label{thm:psu3Embedding}
If $N$ is a compact, orientable, real analytic Riemannian $3$-manifold, then $N$ can be embedded isometrically as a calibrated submanifold in a $9$-manifold $M$ with a $\PSU(3)$-structure with $\rho$ closed.
\end{theorem}
The proof of this theorem will be given in Section~\ref{sec:AbstractEmbedding}. For the moment, we point out that by a real analytic Riemannian manifold we mean a real analytic manifold with a real analytic Riemannian metric.

The only r\^ole of the dimension three in the proof of Theorem~\ref{thm:psu3Embedding}  is that of ensuring that $N$ is parallelizable. In general, one has to take this as a hypothesis. For example, we can consider the analogous problem of $5$-manifolds calibrated by $*\rho$:
\begin{theorem}
\label{thm:psu3Embedding2}
If $N$ is a compact, parallelizable, real analytic Riemannian $5$-manifold, then $N$ can be embedded isometrically as a calibrated submanifold of a $9$-manifold $M$ with a $\PSU(3)$-structure with $*\rho$ closed.
\end{theorem}

\begin{remark}
The form $\rho\in\Lambda^3\su(3)$ can also be defined in terms of the Killing form $B$ of $\su(3)$, as it satisfies (up to a constant)
\[\rho(X,Y,Z)=B([X,Y],Z).\]
More generally, any Lie algebra has a standard $3$-form defined this way, and one can study the associated geometry. However, the results of this paper have no evident application to this general situation, since the form is generally not stable.
\end{remark}
\section{Special geometries, differential graded algebras and exterior differential systems}
\label{sec:Preparatory}
In this section we introduce the exterior differential system associated to special geometries of a specific type, modelled on a certain differential graded algebra (DGA), and we introduce the Cartan-K\"ahler machinery, which will be used in later sections to prove an abstract embedding theorem. As a first application, we determine sufficient conditions for this structure to determine an involutive exterior differential system (implying, in particular, local existence).

We now introduce some notation which will be fixed for the rest of the paper. Let $M$ be  a manifold of dimension $n$, and let $\pi\colon F\to M$ be the bundle of frames. We set $T=\R^n$ as a representation of $\LieG{O}(n)\subset\GL(n,\R)$, so that $TM=F\times_{\GL(n,\R)}T$; in addition, there is a fixed a metric on $T$, which will be freely identified with its dual $T^*$. Let $G$ be a compact Lie group in $\LieG{O}(n)$. The fixed set $(\Lambda^*T)^G$ of $G$ in the exterior algebra over $T$ maps naturally to forms on $F$ by a map
\begin{equation}
\label{eqn:TautologicalMorphism}
\iota\colon (\Lambda^*T)^G\to\Omega(F).
\end{equation}
Notice that $\iota$ does not depend on $G$, but only on $T$ and the tautological form on $F$; we shall sometimes write $\iota_T$ rather than $\iota$ to emphasize the dependence. If $P_G\subset F$ is a $G$-structure, the restriction of the image of $\iota$ to $P_G$ determines a finite dimensional subalgebra of $\Omega(M)$.

Let $\mathcal{I}^G(F)\subset \Omega(F)$ be the ideal generated by $d(\iota((\Lambda^*T)^G))$; by construction, this is a differential ideal, modelled on the graded algebra $(\Lambda^*T)^G$. More generally, let $A$ be a graded subalgebra of $(\Lambda^*T)^G$; we say that a linear endomorphism $f$ of $A$  is a \dfn{differential operator} if $(A,f)$ is a DGA, namely $f^2=0$, the graded Leibnitz rule holds and $f$ is the sum of linear maps
\[f\colon (A\cap \Lambda^pT)\to  (A\cap \Lambda^{p+1})^G.\]
Given the DGA $(A,f)$, one can consider the ideal $\mathcal{I}_{A,f}(F)$ generated by
\[\left\{d\iota(\alpha)-\iota(f(\alpha)), d\iota(f(\alpha))\mid  \alpha\in A\right\}.\]
This is a differential ideal which coincides with $\mathcal{I}^G(F)$ if $(A,f)=((\Lambda^*T)^G,0)$. Since the algebra $A$ is encoded in the definition of $f$, we shall often write $\mathcal{I}_f(F)$ for $\mathcal{I}_{A,f}(F)$.

We fix a torsion-free connection form $\omega$, whose components will be denoted by $\omega_{ij}$, and the tautological form $\theta$ with components $\theta_i$; by the torsion-free assumption,
\begin{equation}
\label{eqn:structure}
d\theta_i=-\sum_{j=1}^n\omega_{ij}\wedge\theta_j.
\end{equation}
Let $Gr_n(TF)$ be the Grassmannian of $n$-planes inside $TF$; by definition, an element of $Gr_n(TF)$ is transverse to $\pi$ if it does not intersect the vertical distribution $\ker\pi_*$. Transverse elements define an open subset $Gr_n(TF,\pi)\subset Gr_n(TF)$, on which we have natural coordinates $(u,p_{ijk})$, where $u$ represents a set of coordinates on $F$ and $p_{ijk}$ are determined by
\begin{equation}
\label{eqn:grassmann}
\omega_{ij}=\sum_{k=1}^n p_{ijk}\theta_k.
\end{equation}
We define $V_n(\mathcal{I}_f)$ as the set of $\pi$-transverse integral elements of dimension $n$, namely those elements $E$ of $Gr_n(F,\pi)$ such that all forms in $\mathcal{I}_f$ restrict to zero on $E$. In terms of the coordinates $(u,p_{ijk})$, we see from \eqref{eqn:structure}, \eqref{eqn:grassmann} that $V_n(\mathcal{I}_f)$ is cut out by affine equations in the $p_{ijk}$ with constant coefficients; hence, it is either empty or a smooth manifold.

The coordinates $p_{ijk}$ determine a natural map
\begin{equation}
 \label{eqn:grassmannian}
\gl(n,\R)\otimes T\to Gr_n(T_uF,\pi).
\end{equation}
Under this map,  $V_n(\mathcal{I}_f)$ corresponds to an affine space
\[Z_f=Z_{A,f}\subset \gl(n,\R)\otimes T.\]
We shall say that an element of $\gl(n,\R)\otimes T$ is an \dfn{integral} of $\mathcal{I}_f$ if it lies in $Z_f$. Notice that the map  \eqref{eqn:grassmannian} depends on the choice of a connection, but not so the space $Z_f$, i.e. whether a point in $\gl(n,\R)\otimes T$ is mapped to an integral.

Define $Z'_f$, $Z''_f$ by the diagram of $G$-equivariant maps
\[
 \xymatrix{
\gl(n,\R)\otimes T\ar[r]^\cong &T\otimes \gl(n,\R)\ar[r]& T\otimes \Lambda^2T\ar[r]& \frac{ T\otimes \Lambda^2T}{\lie{g}\otimes T}\\
Z_f\ar@{^{(}->}[u]\ar[rr]&& Z'_f\ar@{^{(}->}[u]\ar[r] & Z''_f\ar@{^{(}->}[u]}
\]
where the vertical arrows are inclusions and the horizontal arrows are surjective.

We can explain the relation among exterior differential system and special geometries with the following:
\begin{proposition}
\label{prop:SpecialGeometries}
The affine spaces $Z_f$, $Z'_f$, $Z''_f$ are invariant under $G$, and the vector space
\[T\otimes S^2T+ \lie{g}\otimes T\]
is contained in a translate of $Z_f$.
If $E\in Gr_n(T_uF,\pi)$, the following are equivalent:
\begin{itemize}
\item $E$ is an integral element of $\mathcal{I}_f(F)$.
\item Given any connection extending $E$, in the sense that the connection form $\tilde\omega$ satisfies $\ker\tilde\omega_u=E$, the torsion map
\[\Theta\colon P\to T\otimes \Lambda^2T\] satisfies $\Theta(u)\in Z'_f$.
\item Given any local $G$-structure $P_G\subset F$ extending $E$, in the sense that $T_uP_G\supset E$, the intrinsic torsion map
\[\tau\colon P\to \frac{ T\otimes \Lambda^2T}{\lie{g}\otimes T}\]
satisfies $\tau(u)\in Z''_f$.
\end{itemize}
\end{proposition}

In general, one is interested in finding $G$-structures $P_G$ which not only extend an integral element, but are themselves integrals of $\mathcal{I}_f$, i.e. such that $\mathcal{I}_f$ restricts to zero on $P_G$. Then the map $\iota$ descends to an injective DGA morphism $(A,f)\to\Omega(M)$, and Proposition~\ref{prop:SpecialGeometries} tells us that the intrinsic torsion at each point is forced to lie in some space $Z''_f$.

First, we discuss the problem of whether $Z_f$ is empty.
\begin{lemma}
\label{lemma:ZEmpty}
The space $Z_{A,f}$ is non-empty if and only if  $f$ extends to a degree $1$ derivation $\tilde f$ of the graded algebra $\Lambda^*T$. In that case, $Z_{A,f}$ is a translate of $Z_{A,0}$, and $\tilde f$ can be taken to be $G$-equivariant.
\end{lemma}
\begin{proof}
 Take an element $\eta$ of $\gl(T)\otimes T$; its image in $T\otimes\Lambda^2T$ can be wiewed as a linear map $T\to\Lambda^2T$, which extends uniquely to a degree $1$ derivation $d_\eta$ of $\Lambda^*T$ (which fails to be a differential operator since $(d_\eta)^2\neq 0$ in general). It follows from the definition that $\eta$ is in $Z'_{A,f}$ if and only if $d_\eta$  restricted to $A$ coincides with $f$, and two such $\eta$ differ by a derivation  that restricts to zero on $A$, i.e. an element of $Z'_{A,0}$.

For the last part of the statement, observe that the group $\GL(T)$ acts on the space of derivations of degree one on $\Lambda^*T$ by
\[g\cdot f(\alpha)= g(f(g^{-1}\alpha)), \quad \alpha\in \Lambda^*T.\]
By construction, if $f'$ is a derivation that coincides with $f$ on  $(\Lambda^*T)^G$, then any derivation $g\cdot f'$ also restricts to $f$ on $A\subset (\Lambda^*T)^G$. Thus we can define
\[\tilde f=\int_{g\in G} (gf')\mu_g,\]
where $\mu_g$ is the Haar measure, and obtain a $G$-equivariant derivation of degree one which extends $f$.
\end{proof}
In fact, the condition of Lemma~\ref{lemma:ZEmpty} is not trivially satisfied; for an example, we refer to Section~\ref{sec:so3}.

We can now introduce a special class of DGA's for which $Z_f$ has a particularly simple description. Let $\lie{g}^\perp$ denote the orthogonal complement of $\lie{g}$ in $\Lambda^2T$.
\begin{proposition}
\label{prop:cit}
For a given compact Lie group $G$ acting on $T$ and $A=(\Lambda^*T)^G$,  $Z''_0=\{0\}$ if and only if
\[Z_0=(\lie{g}\otimes T)\oplus (T\otimes S^2T).\]
In this case,
\[\codim V_n(\mathcal{I})=\dim T\otimes \lie{g}^\perp.\]
For every differential operator $f$ on $(\Lambda^*T)^G$, either $Z''_f$ is empty or $Z''_f=\{\xi_f\}$, where $\xi_f$ is fixed under the action of $G$.
Moreover  the intrinsic torsion of any $G$-structure $P_G$ is completely determined by the map
\[d\circ\iota\colon (\Lambda^*T)^G\to\Omega(P_G).\]
In particular, a $G$-structure is an integral of $\mathcal{I}_f$ if and only if its intrinsic torsion map
\[P_G\to\lie{g}^\perp\otimes T\]
takes constantly the value $\xi_f$.
\end{proposition}
We say $G$ is \dfn{strongly admissible} if one of the equivalent conditions of Proposition~\ref{prop:cit} hold, and that $P_G$ is a $G$-structure with \dfn{constant intrinsic torsion} if the intrinsic torsion map is constant. Proposition~\ref{prop:cit} suggests that it makes sense to consider $G$-structures with constant intrinsic torsion, for $G$ a strongly admissible group; constant intrinsic torsion geometries are the DGA equivalent of Killing spinors, though more general. Examples of these structures will be given in Section~\ref{sec:SUn} and Section~\ref{sec:oddSUn}.

In the case that $f$ is zero, looking for integral $G$-structures is equivalent to requiring that the holonomy be contained in $G$, as in \cite{Bryant}. However, the same techniques also work in the case of constant intrinsic torsion (that is, when $f$ is non-zero) or in the case of subalgebras $A\subset(\Lambda^*T)^G$. Indeed, in the rest of this section we give a sufficient condition on $(A,f)$ ensuring that $\mathcal{I}_{A,f}$ is involutive.

\bigskip
Now fix an exterior differential system $\mathcal{I}_f$. A subspace $E\subset T_uF$ is integral if $\mathcal{I}_f$ restricts to zero on $E$. If $v_1,\dotsc,v_k$ is a basis of $E$, the space of \dfn{polar equations} of $E$ is
\[\mathcal{E}(E)=\left\{(v_1\wedge\dotsb\wedge v_k)\hook\alpha_u \mid \alpha\in \mathcal{I}_f\right\}.\]
Then the annihilator of $\mathcal{E}(E)$  is the union of all integral elements containing $E$.
Given $E\in V_n(\mathcal{I}_f)$, a flag
\begin{equation}
 \label{eqn:flagE}
\{0\}=E_0\subsetneq\dotsb\subsetneq E_n=E
\end{equation}
is \dfn{ordinary} (resp. \dfn{regular}) if for all $k<n$ (resp. $k\leq n$) the dimension of $\mathcal{E}(E')$ is constant for $E'$ in a neighbourhood of $E_k$ in $Gr_k(F)$. By Cartan's test, the flag \eqref{eqn:flagE} satisfies
\[\codim V_n(\mathcal{I}_f)\geq  c_0 + \dotsb + c_{n-1},\quad c_k=\dim\mathcal{E}(E_k),\]
and equality holds if and only if the flag is ordinary.
Moreover by the Cartan-K\"ahler theorem the terminus $E$ of an ordinary flag is contained in an integral manifold of dimension $n$ (see \cite{BryantEtAl}).

Since the map $\iota$ establishes an isomorphism between $T$ and $E^*$, we can define another flag
\begin{equation}
 \label{eqn:flagW}
\{0\}=W_0\subsetneq\dotsb\subsetneq W_n=T, \quad E_k=\iota(W_k^\perp)^o\cap E
\end{equation}
Here the notation is that if $V\subset W$, we define $V^o\subset W^*$ as the subspace of those elements of $W^*$ that vanish on $V$. We now show that the numbers $c(E_k)$ do not depend on $E$ or $f$, but only on the flag \eqref{eqn:flagW}.
\begin{lemma}
\label{lemma:cw}
Given a compact Lie group $G$ acting orthogonally on $T=\R^n$, and a graded subalgebra $A\subset (\Lambda^*T)^G$, to every $W\subset T$ one can associate a subspace $H(W)\subset\gl(T)$ with the following property. Let $F$ be a frame bundle over an $n$-dimensional manifold $M$, $f$ a differential operator on $A$, $E$ in $V_n(\mathcal{I}_f(F))$, and $E'\subset E$ the subspace corresponding to $W\subset T$. Then the polar space $H(E')$ for the exterior differential system $\mathcal{I}_f(F)$ satisfies
\[H(E')=E'\oplus H(W).\]
\end{lemma}
\begin{proof}
If $E\subset T_uF$, we can think of $\iota$ as a map from $T$ to $T_u^*F$. Since $E\supset E'$ is an integral element, all the forms in $\mathcal{E}(E')$ vanish on $E$. This means that $\mathcal{E}(E')$ does not intersect $\iota(T)$. We can therefore consider the reduced polar equations, namely the image of $\mathcal{E}(E')$ under the map
\[T^*_uF \to \frac{T_u^*F}{\iota(T)}\cong \gl(T)^*.\]
When passing to the quotient, the components $\iota(f(\alpha))$ of the generators vanish, and so the reduced polar equations for $\mathcal{I}_{A,f}(F)$ are the same as those for $\mathcal{I}_{A,0}(F)$. Since \eqref{eqn:structure} does not depend on the point $u$, or the manifold $M$,
their image in $\gl(T)^*$ depends only on $W$. Defining $H(W)$ as the subspace  $\mathcal{E}(E')^o$ of $\gl(T)$, the statement follows.
\end{proof}
Noting that $\dim H(E)+\dim\mathcal{E}(E)=n^2+n$, we set
\[c(W)=n^2+n-\dim H(W).\]
Thus, Lemma~\ref{lemma:cw} asserts that $c(E)$ equals $c(W)$ regardless of $f$, $F$ and $E$.
\begin{definition}
\label{dfn:ordinary}
We say a flag $W_0\subsetneq\dotsb\subsetneq W_n=T$ is \dfn{$\mathcal{I}_A$-ordinary} if
\[c(W_0)+\dotsb+c(W_{n-1})=\dim \frac{\gl\otimes T}{Z_0}.\]
\end{definition}
This is a purely algebraic definition, which does not involve directly an exterior differential system. One can show that the $c(W_i)$ add up to the dimension of $T\otimes\lie{g}^\perp$ if and only if $G$ is strongly admissible and the flag is $\mathcal{I}_A$-ordinary (see Proposition~\ref{prop:SUnIsStronglyAdmissible}).
\begin{proposition}
\label{prop:localexistence}
Given a compact Lie group $G$ acting orthogonally on \mbox{$T=\R^n$}, and a subalgebra $A\subset (\Lambda^*T)^G$, suppose there is an $\mathcal{I}_A$-ordinary flag \[W_0\subsetneq\dotsb\subsetneq W_n=T.\] Then, for every differential operator $f$ on  $A$ which extends to a derivation of degree one of $\Lambda^*T$ and every real analytic $n$-dimensional manifold $M$ with frame bundle $F$, the exterior differential system $\mathcal{I}_f(F)$ is involutive.
\end{proposition}
\begin{proof}
 By Lemma~\ref{lemma:ZEmpty}, $Z_{A,f}$ is a translate of $Z_{A,0}$; in particular,  at each $u\in F$ there is an integral element $E\in V_n(\mathcal{I}_f(F))$. By Lemma~\ref{lemma:cw}, the flag \eqref{eqn:flagW} determines a flag
$E_0\subset\dotsb\subset E_n$
such that
\[c_0+\dotsb+c_{n-1}=\dim \frac{\gl\otimes T}{Z_0}.\]
On the other hand the right-hand side equals the codimension of $V_n(\mathcal{I}_f(F))$, and so by Cartan's test the flag $E_0\subset\dotsb\subset E_n$ is ordinary.
Thus, every integral element is ordinary, and so by the Cartan-K\"ahler theorem it is tangent to a real analytic integral of dimension $n$.
\end{proof}

We can now give a sufficient condition for algebras generated by a single stable form. Recall from \cite{Hitchin:StableForms} that a form $\phi\in\Lambda^pT$ is \dfn{stable} if its infinitesimal orbit $\gl(T)\cdot\phi$ relative to the standard $\GL(T)$ action coincides with $\Lambda^pT$.

More generally, given a subspace $E\subset T$, we shall say that a form $\alpha$ in $\Lambda^pT$ is $E$-stable if for the restriction map $\pi_E\colon \Lambda^pT\to \Lambda^pE$ satisfies
\begin{equation*}
\label{eqn:Estable}
\pi_E(\gl(T)\cdot\phi)=\Lambda^p E.
\end{equation*}
In particular, if $\alpha$ is stable then it is also $E$-stable for all $E\subset T$. The converse is not true: for instance, a $4$-form in $\R^8$ with stabilizer $\Spin(7)$ is not stable, but it is stable for every $E\subsetneq T$, since the restriction to a space of codimension one has stabilizer conjugate to $\Gtwo$.
\begin{proposition}
\label{prop:stable}
Let $\alpha$ be an $E$-stable form on $T$, where $E$ has codimension one in $T$. Then every flag
$E_0\subsetneq\dotsb\subsetneq E_{n}$ in $T$
with $E_{n-1}=E$ is $\mathcal{I}^\alpha$-ordinary.
\end{proposition}
\begin{proof}
Since $\pi_{E_i}$ factors through $\pi_{E_{n-1}}$, $\alpha$ is $E_i$-stable for all $i$. Thus we obtain that $\mathcal{E}(E_i)$ consists of $\binom{i}{p}$ equations. On the other hand, the space $Z^\alpha$ is defined by at most $\binom{n}{p+1}$ equations. By Cartan's inequality,
\[\binom{p}{p}+\binom{p+1}{p}\dotsb+\binom{n-1}{p}\leq \codim Z^{\alpha} \leq \binom{n}{p+1},\]
and a straightforward calculation shows that equality must hold.
\end{proof}
In general, for algebras $A$ generated by a single form, an $\mathcal{I}_A$-ordinary flag may exist even if the form is not $E$-stable for any $E$ of codimension one. For instance, the form $e^{12}+e^{34}$ in $\R^7$ is not $E$-stable for any six-dimensional $E$, but the algebra $A$ it generates admits an $\mathcal{I}_A$-ordinary flag. However, the above proposition applies to some significant cases. Recall from Section~\ref{sec:psu3} that $\R^8$ admits a stable three-form $\phi$ with stabilizer $\PSU(3)$. Indeed  more stable forms exist (see \cite{LePanakVanzura}), but we shall focus on forms with compact stabilizers.
\begin{theorem}
\label{thm:involutive}
If  $A$ is one of
\[(\Lambda^*\R^8)^{Sp(2)Sp(1)}, (\Lambda^*\R^8)^{Spin(7)}, \R\rho\subset\Lambda^*\R^8, \R(*\rho)\subset\Lambda^*\R^8\;\]
then for every real analytic manifold $M$ of the appropriate dimension with frame bundle $F$, the exterior differential system $\mathcal{I}_{A,0}(F)$ is involutive.
Moreover if $A$ is one of
\[(\Lambda^*\R^7)^{\Gtwo},(\Lambda^*\R^{2k})^{SU(k)}, (\Lambda^*\R^{2k+1})^{SU(k)},\]
for every real analytic  manifold $M$ of the appropriate dimension $n$ with frame bundle $F$ and every differential operator $f$ on  $A$ which extends to a derivation of degree one $f$ of $\Lambda^*\R^n$, the exterior differential system $\mathcal{I}_{A,f}(F)$ is involutive.
\end{theorem}
\begin{proof}
Given Proposition~\ref{prop:localexistence}, it suffices to show that a flag satisfying Definition~\ref{dfn:ordinary} exists in each case.
For the first part of the statement, we observe that each algebra is generated by an $E$-stable form and apply Proposition~\ref{prop:stable}.
Concerning the second part, we refer to \cite{Bryant} for the case of $\Gtwo$, and to Lemmas~\ref{lemma:sun},~\ref{lemma:sunodd} for the case of $\SU(n)$.
\end{proof}
\begin{remark}
By the second part of the theorem, it follows that the exterior differential systems associated to nearly-parallel $\Gtwo$-structures, nearly-K\"ahler $\SU(3)$-structures and $\alpha$-Einstein-Sasaki $\SU(k)$-structures are involutive. These structures are also characterized by the existence of a Killing spinor. Thus, this result can be viewed as an ``analogue with torsion'' to Bryant's result in the case of a parallel spinor \cite{Bryant}.
\end{remark}

\begin{remark}
Computing the $c_k$ is in general not a straightforward task. For fixed $G$, one can resort to the use of a computer (see \cite{Conti:SymbolicComputations}). In general,  one has to carry out the computation by hand (See Sections~\ref{sec:SUn}, \ref{sec:oddSUn}).
\end{remark}

\section{The abstract embedding theorem}
\label{sec:AbstractEmbedding}
In this section we formalize the embedding problem introduced in Section~\ref{sec:so3}, and prove an abstract embedding theorem implying Theorem~\ref{thm:so3Embedding}; further applications will be given in later sections.

Let $W\subset T$ be a subspace. The action of $G$ on $T$ restricts to an action of
\[H=G\cap\mathrm{O}(W)\]
on $W$. There is a natural map
\[p\colon (\Lambda^*T)^G\to (\Lambda^*W)^H,\]
which is not in general surjective. Given a differential operator $f$ on $A\subset(\Lambda^*T)^G$, there is at most one differential operator $f_W$ such that the  diagram
\begin{equation}
 \label{eqn:fw}
\begin{aligned}
\xymatrix{
A\ar[d]^p\ar[r]^f& A\ar[d]^p\\
p(A)\ar[r]^{f_W}& p(A)
}
 \end{aligned}
\end{equation}
commutes; the condition for the existence of $f_W$ is that
\[\ker p \subset\ker p\circ f.\]
When $f=0$ the condition is trivially satisfied, and we will denote the induced map by $0_W$. The map $f_W$ can be used to define an induced exterior differential system on certain submanifolds of $N$.

Indeed, let $\iota\colon N\to M$ be a submanifold of $M$, with the same dimension as $W$, and assume the normal bundle is trivial. A principal bundle, whose fibre consists of the group $\GL(T,W)$ of matrices in $\GL(T)$ mapping $W$ to itself, is induced on $N$ by
\[F_N=\{u\in \iota^*F\mid u(W)\subset T_{\pi(u)}N\}.\]
If one fixes a trivialization of the normal bundle, one can view the frame bundle of $N$ as a reduction to $\GL(W)$ of the principal bundle $F_N$; however, it will be more convenient to work with $F_N$. Although $F_N$ is not a $\GL(T,W)$-structure, it carries a tautological form nonetheless, and so it has a tautological morphism of the type \eqref{eqn:TautologicalMorphism}, i.e.
\[
\iota_W\colon  p((\Lambda^*T)^G)\to \Omega(F_N).
\]
Like in Section~\ref{sec:Preparatory}, it follows that $f_W$ determines an exterior differential system on $F_N$ which we denote by $\mathcal{I}_{f_W}(F_N)$.

We wish to relate  integrals of $\mathcal{I}_{f}(F)$ and of $\mathcal{I}_{f_W}(F_N)$. In order to do so, we need to establish a compatibility condition between the embedding $N\subset M$ and the inclusion $W\subset T$.
\begin{definition}
Suppose $N$ is an embedded submanifold of $M$ with the same dimension as $W\subset T$ and $P_G$ is a $G$-structure on $M$; if the intersection
\[\iota^*F_N \cap \iota^* P_G\subset\iota^*F\]
contains a principal bundle $P_H$ on $N$ with fibre $H$, we say that the pair $(N,P_H)$ is embedded in $(M,P_G)$ with type $W\subset T$.
\end{definition}
Notice that in the above definition, the intersection
\[(\iota^*F_N \cap \iota^* P_G)_x=\left\{u\in (P_G)_x\mid u(W)\subset T_{\pi(x)}N\right\}\]
is either empty or a single $\tilde H$-orbit, where
\[\tilde H=G\cap(\mathrm{O}(W)\times\mathrm{O}(W^\perp)).\]
So the first condition required in the definition is that the intersection is never empty, or equivalently, a principal bundle with fibre $\tilde H$. On the other hand, giving a reduction to $H$ of this bundle is the same as giving a trivialization of the normal bundle.

\begin{proposition}
\label{prop:submanifold}
Fix a manifold $M$ and an exterior differential system $\mathcal{I}_f$ on its frame bundle. If a $G$-structure $P_G$ is an integral of $\mathcal{I}_f$ and $(N,P_H)$ is embedded in $(M,P_G)$ with type $W\subset T$, then $\ker p\subset \ker p\circ f$ and $P_H$ is an integral for $\mathcal{I}_{f_W}$.
\end{proposition}
\begin{proof}
Consider the commutative diagram
\[\xymatrix{
A\ar@{^{(}->}[r]\ar[d]^p&(\Lambda^*T)^G\ar[d]^p\ar[r]^{\iota_T} &\Omega(F)\ar[r]\ar[d] &\Omega(P_G)\ar[d]\\
p(A)\ar@{^{(}->}[r]&(\Lambda^*W)^H\ar[r]^{\iota_W} &\Omega(F_N)\ar[r]&\Omega(P_H)
 }
\]
where the straight unlabeled arrows are restriction maps. Let $p(\alpha)$ be a form in $p(A)$. On $P_H$, \[d(\iota_W(p(\alpha)))=d(\iota_T(\alpha))=\iota_T(f(\alpha))=\iota_W(p(f(\alpha))),\]
and so if $p(\alpha)=0$, then $p(f(\alpha))=0$. Since $p\circ f=f_W\circ p$, $P_H$ is an integral of $\mathcal{I}_{f_W}$.
\end{proof}

Conversely, we can pose the following general problem, which we shall call the \dfn{embedding problem}. Let $N$ be a manifold of the same dimension as $W$, and let $F_N$ be the bundle of frames on $N$. Suppose $F_N$ has a $H$-reduction $P_H$ which is an integral of $\mathcal{I}_{f_W}(F_N)$. Can we embed $(N,P_H)$ into a pair $(M,P_G)$ with type $W\subset T$, in such a way that $P_G$ is an integral of $\mathcal{I}_{f}$?

\begin{remark}
In principle one could also consider submanifolds $N$ with non-trivial normal bundle. So far, this would amount to replacing the  structure group $H$ with $\tilde H$.
However, the main result of this section makes use of Cartan-K\"ahler theory, which requires enlarging one dimension at a time. Thus, the manifolds one obtains have trivial normal bundle and the $\tilde H$-structure induced by the embedding reduces to $H$. For this reason, the theorem will not apply to $\tilde H$-structures which do not admit a $H$-reduction. See also the comments in \cite{Bryant:Calibrated} and the remark at the end of this section.
\end{remark}

It was shown in \cite{Bryant:NonEmbedding} that one cannot solve the embedding problem in general in the non-real-analytic setting. However, keeping the real analytic assumption, we will prove an abstract embedding theorem, which generalizes the results mentioned in Section~\ref{sec:so3}, does not reference any specific instance of $G$ and $T$, and holds in the constant intrinsic torsion case as well.

\begin{lemma}
\label{lemma:extends}
Assume $Z_f$ is not empty.
For ever choice of  $W\subset T$, the linear projection $\gl(T)\otimes T\to\gl(W)\otimes W$ induces by restriction a surjective map \mbox{$Z_f\to Z_{f_W}$}.
\end{lemma}
\begin{proof}
It follows from the diagram \eqref{eqn:fw} that the image of $Z_f$ is contained in $Z_{f_W}$; in particular $Z_{f_W}$ is not empty. Since the projection is linear, the image of $Z_f$ has the same dimension as the image of $Z_0$. Moreover, as $Z_{f_W}$ is not empty, its dimension does not depend on $f$. Summing up, it suffices to prove the statement for $f=0$.

The rest of the proof, while algebraic in nature, is best explained working at a point $u$ of the frame bundle $F$.
Let $W$ have dimension $k$ and $T$ have dimension $n$.  Consider the natural map
\[\gl(W)\otimes W\to \gl(T)\otimes T\to Gr_n(T_uF,\pi)\xrightarrow{q} Gr_k(T_uF,\pi),\]
where the map $q$ is given by the choice of $W\subset T$ and the correspondence of \eqref{eqn:flagW}. We claim that the image $E$ of an element of $Z_{0_W}$ under this map is an integral element of $\mathcal{I}_{A,0}(F)$.
Indeed, let $\alpha$ be in $A$; we can decompose it according to $T=W\oplus W^\perp$, so that
\[\alpha=p(\alpha)+\sum \beta_i\wedge\gamma_i, \quad \beta_i\in\Lambda^*W^\perp, \gamma_i\in\Lambda^*W.\]
Restricting to the space $E$,  $d\iota(p(\alpha))$ is zero because $E$ is the image of an element in $Z_{0_W}$, $\iota(\beta_i)$ is zero  because of how $q$ is defined, and finally $d\iota(\beta_i)$ is zero because $E$ is the image of an element in $\gl(W)\otimes W$. Thus, $d\iota(\alpha)$ is zero on $E$ for all $\alpha$, and $\pi$-transverse integral elements of $\mathcal{I}_{0_W}$ are mapped to $\pi$-transverse integral elements.

Therefore, we obtain a commutative diagram
\[\xymatrix{
Z_0\ar[d]\ar[r] & V_n(\mathcal{I}_{A,0}(F))\ar[d]^q\\
Z_{0_W}\ar[r] & V_k(\mathcal{I}_{A,0}(F))
 }
\]
where the bottom arrow is the map constructed above; it suffices to show that the map $q$ is surjective. In other words, we must prove that for the exterior differential system $\mathcal{I}_{A,0}(F)$, every $\pi$-transverse integral element $E\subset T_uF$ of dimension $k<n$  is contained in a $\pi$-transverse integral element of dimension $n$.

Let $e$ be a generator of  $\Lambda^k E$, and set
\[\mathcal{\tilde E}(E)=\{e\hook \phi_u\mid \phi\in\mathcal{I}^G(F)\}\;.\]
By construction, $\mathcal{\tilde E}(E)$ is closed under wedging by forms in the exterior algebra $\Lambda E^o$ over the space $E^o\subset T^*_uF$ of $1$-forms that vanish on $E$.  The polar equations $\mathcal{E}(E)$ are the subspace of $\mathcal{\tilde E}(E)$ consisting of forms of degree one; a vector space $E'\supset E$ of dimension $k+1$ is integral if and only if all forms in $\mathcal{E}(E)$ vanish on $E'$. Introducing the space of ``horizontal'' one-forms
\[L=(\ker\pi_{*u})^o,\]
we now prove that $\mathcal{E}(E)$ is transverse to $L$. Indeed, let $\Theta$ be a generator of $\Lambda^nL$, and let $\Theta_o$ be a generator of $\Lambda^{n-k}(E^o\cap L)$. Then $\Theta_o$ does not lie in $\mathcal{\tilde E}(E)$, because otherwise $\mathcal{I}_{A,0}$ would contain an element of the form $\Theta+\beta$, $e\hook\beta=0$; this is absurd, because it implies that the zero element in $\gl(T)\otimes T$ is not an integral of $\mathcal{I}_{A,0}$. Thus $\mathcal{\tilde E}(E)$, which  is closed under wedging with forms in $\Lambda E^o$, does not contain any $1$\nobreakdash-form in $E^o\cap L$. On the other hand $E$ is integral, so  $\mathcal{E}(E)$ is contained in $E^o$ and
\[\mathcal{E}(E)\cap L=\mathcal{E}(E)\cap (E^o\cap L)=\{0\}\;.\]
Thus, a basis of $L$ restricts to linearly independent elements on $\mathcal{E}(E)^o$; this implies the existence of a $\pi$-transverse integral $E'\supset E$ of dimension $k+1$, and by repeating the argument we obtain the required integral element of dimension $n$.
\end{proof}

In Proposition~\ref{prop:localexistence} we introduced a sufficient condition for $\mathcal{I}_f$ to be involutive. An involutive exterior differential system is a necessary ingredient for our embedding theorem, but the condition needs to be refined slightly.  We say that $W\subset T$ is \emph{relatively admissible} if
there is an $\mathcal{I}^G$-ordinary flag
\begin{equation}
 \label{eqn:relativelyadmissible}
W_0\subsetneq\dotsb\subsetneq W_{n}=T, \quad W_k=W.
\end{equation}

We can now prove the main result of this section, which also gives Theorems~\ref{thm:psu3Embedding},~\ref{thm:psu3Embedding2}  as corollaries.
\begin{theorem}[Abstract Embedding Theorem]
\label{thm:AbstractEmbedding}
Let $G\subset \LieG{O}(T)$ be a compact group, and let $W\subset T$ be a relatively admissible subspace, with $H=G\cap \LieG{O}(W)$. Let $f$ be a differential operator on $A\subset (\Lambda^*T)^G$ which extends to a derivation of degree $1$ on $\Lambda^*T$ and assume $\ker p\subset \ker f\circ p$; let $f_W$ be the differential operator induced by the diagram~\eqref{eqn:fw}. Then every real analytic pair $(N,P_H)$, with $N$ a manifold of the same dimension as $W$ and $P_H$  an integral of $\mathcal{I}_{p(a),f_W}$,  can be embedded in a pair $(M,P_G)$ with type $W\subset T$, where $P_G$ is an integral of $\mathcal{I}_{A,f}$.
\end{theorem}

\begin{proof}
We shall prove that if $W'=W_{k+1}$, with $W_{k+1}$ defined by the flag \eqref{eqn:relativelyadmissible}, then $(N,P_H)$ can be embedded with type $W\subset W'$ in a real analytic pair $(N',P_K)$, where
\[K=G\cap\mathrm{O}(W'),\]
and $P_K$ is an integral of $\mathcal{I}_{f_{W'}}$. By definition, either $W'=T$ and there is nothing left to prove, or $W'\subset T$ is relatively admissible and $(N',P_K)$ satisfies the hypotheses of the theorem. Therefore, we can repeat the argument, obtaining a chain of embedded pairs. Now, if $(N,P_H)$ is embedded in $(N',P_{H'})$ with type $W\subset W'$, and $(N',P_{H'})$ is embedded in $(M,P_G)$ with type $W'\subset T$, then $(N,P_H)$ is embedded in $(M,P_G)$ with type $W\subset T$. Therefore, the statement will follow.

Having reduced the problem to one-dimensional steps, we proceed by adapting the proof of \cite{ContiSalamon}. Suppose $P_H$ is an integral of $\mathcal{I}_{f_W}$ on $F_N$. Let $M=N\times \R^{n-k}$; denote by $q\colon M\to N$ the natural projection. The frame bundle $F_N$ of $N$ pulls back to a $\GL(W)$-structure $q^*F_N$  on $M$. We claim that $(q^*P_H)|_{M\times\{0\}}$ is an integral of the exterior differential system $\mathcal{I}_{f}(q^*F_N)$. Indeed, consider the commutative diagram
\[\xymatrix{
A\ar[r]^{\iota_{T}}\ar[d]^{p} &\Omega\left(q^*F_N|_{M\times\{0\}}\right)\ar[r]&\Omega\left(q^*P_H|_{M\times\{0\}}\right)\\
p(A)\ar[r]^{\iota_W} &\Omega(F_N)\ar[r]\ar[u]_\cong&\Omega(P_H)\ar[u]_\cong
 }
\]
where the unlabeled arrows are restriction maps.
For any $\alpha\in (\Lambda^*T)^G$, we see that the form upstairs
\[d\iota_{T}(\alpha)-\iota_{T} f(\alpha)\]
can be identified with the form downstairs
\[d\iota_{W}(p(\alpha))-\iota_{W} f_{W}(p(\alpha)),\]
which vanishes on $P_H$. Therefore $(q^*P_H)|_{M\times\{0\}}$ is an integral of $\mathcal{I}_{f}(q^*F_N)$, and in particular an integral of $\mathcal{I}_{f}(q^*F_M)$, where $F_M$ denotes the frame bundle of $M$.

Now fix a torsion-free connection on $F_N$,  pull it back to $q^*F_N$ and extend it to $F_M$. At a point $u$ of $P_H$, we can write
\[T_uF_N=W\oplus \gl(W).\]
Choose an $E\subset T_uF_N$ such that
\[T_uP_H=E\oplus \lie{h}.\]
Then the usual map \eqref{eqn:grassmannian} identifies $E$ with an element of $Z_{f_W}\subset\gl(W)\otimes W$, which we also denote by $E$. Passing to the diffeomorphic bundle  $q^*P_H|_{M\times\{0\}}$, the point $u$ corresponds to a point $(u,0)$, and the splitting
\[T_{(u,0)}q^*F_N=T\oplus \gl(W)\]
is compatible with the splitting of $T_uF_N$; thus, the tangent space to $(q^*P_H)|_{M\times\{0\}}$ at ${(u,0)}$ can also be identified with $E+\lie{h}$, where $E$ is in $Z_{f_W}$. By Lemma~\ref{lemma:extends}, $E$ is the image of an element $\tilde E$ of $Z_f$ under the map
$\gl(T)\otimes T\to\gl(W)\otimes W$. In terms of bundles, this means that
\[T_{(u,0)}q^*P_H\subset \tilde E\oplus \lie{h}\]
where $\tilde E$ is in $V_n(\mathcal{I}_f(q^*F_N))$. Since $W\subset T$ is relatively admissible, and using Lemma~\ref{lemma:cw}, we can find a  flag
\[E_0\subset \dotsb\subset E_{n}=\tilde E, \quad E=E_k;\]
which satisfies Cartan's test, and is therefore ordinary. In particular, $E$ is regular.

By the $G$-invariance (and therefore $H$-invariance) of $\mathcal{I}_f$,  a $\pi$-transverse element $E_0$ is integral if and only if $E_0+\lie{h}$ is integral, and the polar spaces are related by
\[H(E_0)\oplus\lie{h}=H(E_0+\lie{h}).\]
In order to apply the Cartan-K\"ahler theory, we must quotient out the $H$ invariance. Thus, we consider the quotient $F_M/H$. Then the  manifold $q^*P_H/H\mid_{M\times\{0\}}$ is an integral of $\mathcal{I}_f(F_M/H)$. Moreover since $\lie{h}$ is contained in any polar space, the ordinary flag
$E_0\subsetneq \dotsb\subsetneq E_n$ determines an ordinary flag $E_0'\subsetneq \dotsb \subsetneq E_n'$ on the quotient $F_M/H$.
In particular, $E_k'$ is regular, and since this holds at all $u$, the integral manifold
$(q^*P_H/H)|_{N\times\{0\}}$
is regular.

The Cartan-K\"ahler theorem requires a ``restraining manifold'', namely a real analytic submanifold $R$ of $F_M/H$ which intersects each polar space of $(q^*P_H/H)|_{N\times\{0\}}$ transversely, and such that the codimension of the former is the extension rank of the latter. We shall define $R$ as the quotient of a $H$-invariant submanifold $\tilde R$ of $F_M$ with the same properties. The extension rank can be defined as
\[r(E+\lie{h})=\dim H(E+\lie{h})-\dim (E+\lie{h})-1,\]
and since $H(E)\supset\lie{h}\oplus E_n$, it satisfies
\[r(E+\lie{h})=\dim H(E)-k-1-\dim\lie{h} \geq n-k-1\geq 0.\]
By Lemma~\ref{lemma:cw}, the polar space of $E$ satisfies
$H(E)=\tilde E\oplus H(W)$. Let $V$ be a $H$-invariant complement of $H(W)$ in $\gl(T)$, and let $U$ be the image under the exponential map of a small neighbourhood of zero in $V$, so that $U$ is invariant under the adjoint action of $H$ and the product map
$U\times H\to \GL(T)$  is an embedding.

Denote by  $N'$ the subset
\[ \{(x,y_1,\dotsc, y_{n-k}\in M\times\R^{n-k}, \mid y_2=\dotsc=y_{n-k}=0\}.\]
Then the manifold
\[\tilde R=(q^*P_H\cdot U)|_{N'}\]
is a fibre bundle over $N'$ with fibre $H\times U$. However by construction $H\times U=U\times H$, and so $\tilde R$ is  $H$-invariant; its  codimension is
\[\dim V+n-k-1- \dim\lie{h}=\dim H(E)-k-1-\dim\lie{h}=r(E+\lie{h}).\]
Moreover, at a point $(u,0)$ of $q^*P_H$ we have
\[T_{(u,0)}\tilde R=T_{(u,0)}q^*P_H|_{N\times\{0\}} \oplus \left\langle\frac\partial{\partial y_1}\right\rangle \oplus V,\]
so
\[T_{(u,0)}\tilde R+H(E)=E+H(W)+V=T_{(u,0)}F_M.\]
Then the quotient $R=\tilde R/H$ is the restraining manifold we need.

Applying the Cartan-K\"ahler theorem to $F/H$, and taking the preimage in $F$, we obtain a $H$-invariant integral manifold $Q$, $(q^*P_H)|_{N\times\{0\}}\subset Q\subset \tilde R$ of dimension $k+1+\dim \lie{h}$. By construction,
\[T_{(u,0)}Q=E_{k+1}\oplus \lie{h}\]
where $E_{k+1}$ is $\pi$-transverse. Up to restricting $N'$,  $P_K=Q\cdot K$ is a $K$-structure on $N'$ as well as an integral manifold of $\mathcal{I}_f$, and therefore of $\mathcal{I}_{f_{W'}}$; moreover,  $(N,P_H)$ is embedded in $(N',P_{H'})$ with type $W\subset W'$.
\end{proof}
As a corollary, we obtain the following:
\begin{proof}[Proof of Theorem~\ref{thm:psu3Embedding}]
Let $A$ be the algebra generated by $\rho$, and let $W\subset\R^8$ be a subspace calibrated by $\rho$, i.e. a maximum point for \eqref{eqn:comass}. Since $\rho$ is stable, $W\subset\R^8$ is relatively admissible with respect to $\mathcal{I}_A$. Being three-dimensional, the manifold $N$ is parallelizable, and compactness ensures that a real analytic parallelism exists (see \cite{Bryant:Calibrated,Bochner}), which can be taken to be orthonormal. Thus, $N$ has a real analytic $\{e\}$-structure, which is an integral of $\mathcal{I}_{0_W}$ because $A$ is generated by a three-form.
By Theorem~\ref{thm:AbstractEmbedding} the manifold $N$ can be embedded with type $W\subset\R^8$ in a manifold $M$ with a $\PSU(3)$-structure whose associated form is closed. Since $W$ is a calibrated space, $N$ is a calibrated submanifold.
\end{proof}
\begin{proof}[Proof of Theorem~\ref{thm:psu3Embedding2}]
Analogous to the proof of Theorem~\ref{thm:psu3Embedding}, except that parallelizability is now part of the hypothesis.
\end{proof}

Notice that Theorem~\ref{thm:AbstractEmbedding} applies to all the geometries appearing in Theorem~\ref{thm:involutive}, and in particular to
\begin{equation}
\label{eqn:ParallelSpinor}
(\Lambda^*\R^{2n})^{\SU(n)}, (\Lambda^*\R^7)^\Gtwo, (\Lambda^*\R^8)^{\Spin(7)}.
\end{equation}
In these cases every codimension one $W\subset T$ is relatively admissible, because the structure group acts transitively on the sphere.
So, in some sense Theorem~\ref{thm:AbstractEmbedding} answers the generalized Killing spinor embedding problem mentioned in the introduction. In fact, by \cite{Wang:ParallelSpinorsAndParallelForms} the holonomy of a simply-connected, irreducible Riemannian manifold with a parallel spinor is either $\Gtwo$, $\Spin(7)$, $\SU(n)$ or $\LieG{Sp}(n)$. Since $\LieG{Sp}(n)$ is contained in $\SU(2n)$, a real analytic Riemannian manifold with a generalized Killing spinor $\psi$ can be isometrically embedded in an irreducible Riemannian manifold with a parallel spinor restricting to $\psi$ if and only if the $G$-structure defined by $\psi$ is an integral of $\mathcal{I}_{p(A),0_W}$, where $A$ is a DGA in \eqref{eqn:ParallelSpinor}. Other examples will be given in Sections \ref{sec:hypo}, \ref{sec:go}.
\begin{remark}
Theorem~\ref{thm:AbstractEmbedding} does not require that $G$ be admissible or strongly admissible. However, if the group is strongly admissible Proposition~\ref{prop:cit} gives  a much stronger interpretation in terms of intrinsic torsion.
\end{remark}

\begin{remark}
One can weaken the assumptions slightly and consider submanifolds whose normal bundle is flat rather than trivial; in other words, one replaces the group $H$ with a larger group $H'\subset G$ that acts discretely on $W^\perp$. The idea is that any $H'$-structure on $N$ reduces to $H$ when pulled back to the universal covering $\tilde N$. One can then apply the theorem to $\tilde N\subset \tilde M$, and obtain the embedding $N\subset M$ by means of a quotient. Indeed, with notation as in the proof of the theorem, $\pi_1(N)$ acts properly discontinuosly on $q^*P$ and therefore on $R$. This action is compatible with the product action on $N\times\R$. Under suitable assumptions on $N$, e.g. if $\pi_1(N)$ is finite or $N$ is compact, we can restrict $P_{H'}$ and $N'$ to make them $\pi_1(N)$-invariant; moreover, the action is properly discontinuos because it is on the base.
\end{remark}

\section{The structure group $\SU(n)\subset \LieG{O}(2n)$}
\label{sec:SUn}
In this paper we have considered $\SU(n)$ as a subgroup of $\LieG{O}(2n)$ or as a subgroup of $\LieG{O}(2n+1)$; either way it is strongly admissible. In this section we consider it as a subgroup of $\LieG{O}(2n)$, and show that it admits an $\mathcal{I}^{\SU(n)}$-ordinary flag, so that Theorem~\ref{thm:AbstractEmbedding} applies. This result will be used in Section \ref{sec:hypo}.

Let $T=\R^{2n}$,  $n\geq 2$, with basis $\theta_1,\dotsc,\theta_{2n}$, and identify $\SU(n)$ with the subgroup of $\LieG{O}(T)$ which fixes the real two-form and complex $n$-form
\begin{equation}
\label{eqn:SUnForms}
F=\theta^{12}+\dotsc + \theta^{2n-1,2n}, \quad \Omega=\Omega^+ + i\Omega^-=(\theta_1+i\theta_2)\dotsm (\theta_{2n-1}+i\theta_{2n}).
\end{equation}

Define a flag $E_0\subsetneq\dotsb\subsetneq E_{2n}=T$, where for $0\leq k\leq n$
\[E_k=\Span{\theta_1,\dotsc,\theta_{2k-1}},\quad E_{n+k}=E_n\oplus\Span{\theta_2,\dotsc, \theta_{2k}}.\]
Let $S=\gl(T)$, with basis $\omega_{ij}$. If $F$ is the frame bundle over a $2n$-dimensional manifold with a fixed torsion-free connection, the direct sum $T\oplus S\subset\Omega(P)$ maps isomorphically to $T^*F_u$ for any $u\in F$. In addition, exterior differentiation on $F$ restricts to an operator
\begin{equation}
 \label{eqn:formald}
d\colon \Lambda T\to S\otimes\Lambda T,\quad d(\theta_i)=\sum_{1\leq j\leq 2n}\omega_{ij}\theta_j.
\end{equation}
In order to compute the numbers $c(W_k)$ defined after Lemma~\ref{lemma:cw} we can work in the exterior algebra over $T\oplus S$, using the formal operator $d$ defined by \eqref{eqn:formald}.
\begin{lemma}
\label{lemma:sun}
With respect to the exterior differential system $\mathcal{I}^{\SU(n)}$, the flag defined above satisfies
\[c(W_k)=\begin{cases}
	\binom{k}{2}& k<n\\
	2-3n^2+4nk-n-\binom{k}2& n\leq  k<2n\\
	1+3n+6\binom{n}2&k=2n\\
\end{cases}.\]
In particular, \[\sum_{k=0}^{2n-1} c_k = 2n(n^2-n+1)=\dim \R^{2n}\otimes\frac{\so(2n)}{\su(n)}.\]
\end{lemma}
\begin{proof}
Let $E\subset T$ be the vector space spanned by some subset of $\{\theta_1,\dotsc,\theta_{2n}\}$. We write
\[\Lambda^n(E)=L(E)\oplus C(E)\]
where $L(E)$ is spanned by elements $\theta_1\wedge\dotsb \wedge \theta_n$ such that the set $\{\theta_i\}$ does not contain any pair $\{\theta_{2j-1},\theta_{2j}\}$, and $C(E)$ is spanned by elements $\theta_1\wedge\dotsb \wedge \theta_n$ such that the set $\{\theta_i\}$ contains at least one pair $\{\theta_{2j-1},\theta_{2j}\}$. Similarly, we write
\[\Lambda^2(E)=L^2(E)\oplus C^2(E),\]
where $C^2(E)$ is spanned by elements of the form $\theta_{2j-1}\wedge \theta_{2j}$ and $L^2(E)$ is its natural complement.
In addition, we further decompose $L(E)$ as
\[L(E)=L^+(E)\oplus L^-(E)\]
where  $\theta_{i_1}\wedge\dotsb\wedge \theta_{i_n}$ lies in $L^+(E)$ or $L^-(E)$  according to whether the number of odd $i_k$ is even or odd.

Since we have fixed a basis, the space $T\oplus S$ can be identified with its dual; thus, we have natural maps
\[\phi\colon \Lambda^2 T\to S,\quad \phi(v\wedge w)=dF(v,w),\]
\[\phi^\pm\colon \Lambda^n T\to S,\quad \phi(v_1\wedge\dotsm \wedge v_n)=d\Omega^\pm(v_1,\dotsc,v_n),\]
where as usual $\Omega^\pm$ are real forms such that $\Omega=\Omega^++i\Omega^-$. The space of polar equations of $E$ can be expressed as
\[\mathcal{E}(E)=\phi(\Lambda^2(E))+ \phi^+(\Lambda^n(E))+ \phi^-(\Lambda^n(E)).\]
To compute its dimension, we decompose $S$ as
\[S=S_1\oplus S_2\oplus S_3\]
where
\[S_1=\Span{\omega_{ii}\mid 1\leq i\leq 2n},\quad S_2=\Span{\omega_{2k,2k-1},\omega_{2k-1,2k}\mid 1\leq k\leq n},\]
and $S_3$ is spanned by the $\omega_{ij}$ that do not lie in $S_1\oplus S_2$.

As a first observation, notice that $\phi$ can be represented as matrix multiplication by the invertible square matrix representing the complex structure on $\R^{2n}$, and so it is injective. With respect to the above decomposition,
\begin{equation}
\label{eqn:imagesplits1}
\phi(L^2(E))\subseteq S_3, \quad \phi(C^2(E))\subseteq  S_1,
\end{equation}
and more explicitly
\begin{equation}
 \label{eqn:explicitly}
\phi(C^2(E_k))=\Span{\omega_{2h-1,2h-1}+\omega_{2h,2h}}_{1\leq h\leq k-n}.
\end{equation}
Similarly, we have
\begin{equation}
\label{eqn:imagesplits2}
\phi^\pm(L^\pm(E))\subseteq S_1,\quad
\phi^\pm(L^\mp(E))\subseteq S_2,
 \quad  \phi^\pm(C(E))\subseteq S_3.
\end{equation}
Next observe that $\phi^\pm$ is injective on $L^\pm(E)$. Now write
\[\Omega=(\theta_1+i\theta_2)\wedge\Omega_{n-1},\quad \Omega_{n-1}=(\theta_3+i\theta_4)\dotsm (\theta_{2n-1}+i\theta_{2n}).\]
If $\theta_1\wedge\alpha$ is in $L$, where $\alpha$ has no component in $\theta_1\wedge\Lambda T$,
then
\[(\theta_1\wedge\alpha)\hook d\Omega +i (\theta_2\wedge\alpha)\hook d\Omega = -(\omega_{11}-\omega_{22}+i\omega_{21}+i\omega_{12}) \alpha\hook \Omega_{n-1}.\]
So if  $\theta_1\wedge\alpha$ is in $L^+$ then $\alpha\hook\Omega_{n-1}$ is a real non-zero number, which we can normalize to $1$, finding that
\begin{align*}
\phi^+(\theta_1\wedge\alpha) &= \phi^-(\theta_2\wedge\alpha) -\omega_{11}+\omega_{22}\\
\phi^-(\theta_1\wedge\alpha) &= -\phi^+(\theta_2\wedge\alpha) -(\omega_{21}+\omega_{12})
\end{align*}
If $E=E_k$, $k>n$ we can compute $\mathcal{E}(E_k)$ in three steps.

(\emph{i}) The space $\mathcal{E}(E_k)\cap S_1$ contains $\omega_{11}-\omega_{22}$, and therefore, by \eqref{eqn:explicitly} it contains both $\omega_{11}$ and $\omega_{22}$. By relabeling the indices, it follows that $\mathcal{E}(E_k)\cap S_1$ contains the elements
\[ \omega_{2h-1,2h-1},\omega_{2h,2h},\quad 1\leq h\leq k-n.\]
 In addition $\mathcal{E}(E)$ contains
\[\omega_{11}+\omega_{33}+\dotsc + \omega_{2n-1,2n-1}=-\phi^+(\theta_1\wedge\dotsb\wedge \theta_{2n-1}).\]
Thus
\[\dim \mathcal{E}(E_k)\cap S_1=\begin{cases}
	0& k<n\\
	2(k-n)+1& n\leq k<2n\\
	2n&k=2n
                                \end{cases}\]

(\emph{ii}) The space $\mathcal{E}(E_k)\cap S_2$ contains $\omega_{21}+\omega_{12}$,
and
\[\mathcal{E}(E_k)\cap S_2=\Span{\omega_{21}+\omega_{12}}+\phi^+(L_-(E)).\]
By relabeling the indices, it follows that $\mathcal{E}_n(E_k)\cap S_2$ contains the element
\[ \omega_{2h-1,2h}+\omega_{2h,2h-1},\quad 1\leq h\leq k-n.\]
 In addition $\mathcal{E}(E_k)$ contains
\[\omega_{21}+\omega_{43}+\dotsc + \omega_{2n,2n-1}=-\phi^-(\theta_1\wedge \theta_3\wedge \dotsb\wedge \theta_{2n-1}).\]
Thus
\[\dim \mathcal{E}(E_k)\cap S_2=\begin{cases}
	0& k< n\\
	(k-n)+1& n\leq k\leq 2n
                                \end{cases}\]

(\emph{iii}) Finally, $\phi^\pm_n(C(E))$ is spanned by the elements $\phi^\pm(\alpha_{ij})$, where
\[\alpha_{ij}=\theta_{2i-1}\hook \left(\theta_{2j}\wedge \theta_1\wedge \dotsb\wedge \theta_{2n-1}\right), \quad i\neq j.\]
For $1\leq i, j\leq n$, define
\[S_{ij}=\Span{\omega_{hk}\mid \left\{\left[\frac{h+1}2\right],\left[\frac{k+1}2\right]\right\}=\{i,j\}}.\]
Then
\[S_3=\bigoplus_{i<j} S_{ij},\]
and
$\phi^\pm(\alpha_{ij})$ is in $S_{ij}$; more precisely
\begin{align*}
\phi^+(\alpha_{ij})&=\omega_{2i-1,2j}-\omega_{2i,2j-1},& \phi^-(\alpha_{ij})&=\omega_{2i,2j}-\omega_{2i-1,2j-1},\\
\phi^+(\alpha_{ji})&=\omega_{2j-1,2i}-\omega_{2j,2i-1},& \phi^-(\alpha_{ji})&=\omega_{2j,2i}-\omega_{2j-1,2i-1}.
\end{align*}
On the other hand
\begin{align*}
\phi(\theta_{2i-1}\wedge\theta_{2j-1})&=\omega_{2i,2j-1}-\omega_{2j,2i-1},&
\phi(\theta_{2i-1}\wedge\theta_{2j})&=\omega_{2i,2j}+\omega_{2j-1,2i-1},\\
\phi(\theta_{2i}\wedge\theta_{2j-1})&=-\omega_{2i-1,2j-1}-\omega_{2j,2i},&
\phi(\theta_{2i}\wedge\theta_{2j})&=-\omega_{2i-1,2j}+\omega_{2j-1,2i}.
\end{align*}
These elements span a six-dimensional space, so $S_{ij}\cap\mathcal{E}(T)$ has dimension $6$. In general,
the dimension of $\mathcal{E}(E_k)\cap S_{ij}$ only depends on how many of the elements $\theta_{2i-1}$, $\theta_{2i}$, $\theta_{2j-1}$, $\theta_{2j}$ lie in $E_k$.
It follows that, if $1\leq i<j\leq n$:
\[\dim\mathcal{E}(E_k)\cap S_{ij}=\begin{cases}
0&k<j\\
1& j\leq k\leq n\\
1&0\leq k-n<i\\
4&i\leq k-n<j\\
6& i<j\leq k-n \end{cases}
\]
Therefore
\[\dim \mathcal{E}(E_k)\cap S_3=\begin{cases}
	\binom{k}{2}& k\leq n\\
	6\binom{k-n}2+4(k-n)(2n-k)+\binom{2n-k}2& n\leq k\leq 2n\\
                                \end{cases}\]
Since by \eqref{eqn:imagesplits1} and \eqref{eqn:imagesplits2}
\[\mathcal{E}(E_k)=(\mathcal{E}(E_k)\cap S_1) \oplus (\mathcal{E}(E_k)\cap S_2)\oplus (\mathcal{E}(E_k)\cap S_3),\]
the statement follows.
\end{proof}

\section{The structure group $\SU(n)\subset \LieG{O}(2n+1)$}
\label{sec:oddSUn}
In this section we consider $\SU(n)$ as a structure group in dimension $2n+1$. We show that it is strongly admissible, and we construct an $\mathcal{I}^{\SU(n)}$-ordinary flag, so that Theorem~\ref{thm:AbstractEmbedding} applies. This result will be used in Section~\ref{sec:go}.

Let $T=\R^{2n+1}$; the space $(\Lambda^*T)^{\SU(n)}$ is spanned by $\alpha$, $F$ and $\Omega^\pm$, where
\begin{equation}
\label{eqn:oddSUnForms}
\begin{gathered}
\alpha=\theta^{2n+1},\quad F=\theta^{12}+\dotsb+\theta^{2n-1,2n},\\
\Omega=\Omega^+ + i\Omega^-=(\theta^1+i\theta^2)\dotsm (\theta^{2n-1}+i\theta^{2n}).
\end{gathered}
\end{equation}
We start by constructing an $\mathcal{I}^{\SU(n)}$-ordinary flag. Let $E_0\subsetneq\dotsb\subsetneq E_{2n+1}$ be the flag such that, for $0\leq k\leq n$,
\begin{gather*}
 E_k=\Span{\theta_1,\dotsc,\theta_{2k-1}}, \quad 0\leq k\leq n\\
E_{n+k}=E_n\oplus\Span{\theta_2,\dotsc, \theta_{2k}}, \quad 1\leq k<n\\
 E_{2n}=E_{2n-1}\oplus\Span{\theta_{2n+1}}, \quad E_{2n+1}=T
\end{gather*}

\begin{lemma}
\label{lemma:sunodd}
With respect to $\mathcal{I}^{\SU(n)}$, the flag $E_0\subsetneq\dotsb\subsetneq E_{2n+1}$ satisfies
\[c(E_k)=\begin{cases}
	\binom{k}{2}+k& k<n\\
	2+\binom{n}2+n& k=n\\
	2-3n^2+4nk-n-\binom{k}2+k& n< k<2n\\
	3n+2-3n^2+3n(2n-1)&k=2n\\
\end{cases}.\]
In particular the sum $c_0+\dotsb+c_{2n}$ equals $2n^3+3n^2+3n+1$.
\end{lemma}
\begin{proof}
 We proceed like in the proof of Lemma~\ref{lemma:sun}. We have to consider however an extra map
\[\phi^\alpha\colon \Lambda^1(E)\to S, \phi^\alpha(v)=v\hook d\theta_{2n+1}.\]
For fixed $k$ with $k<2n$, we define $\mathcal{\tilde E}(E_k)$ as the space of polar equations with respect to $\mathcal{I}^{\SU(n)}$, and view the space of polar equations computed in Lemma~\ref{lemma:sun} as a subspace $\mathcal{E}(E_k)\subset \mathcal{\tilde E}(E_k)$; we obtain
\[\mathcal{\tilde E}(E_k)=\mathcal{E}(E_k)\oplus\Span{\omega_{2n+1,i}\mid \theta_i\in E_k}.\]
This accounts for the values of $c_k$ for $k<2n$. For $k=2n$, we can see directly that
\[\mathcal{\tilde E}(E_{2n})=\mathcal{\tilde E}(E_{2n-1})\oplus \Span{\omega_{i,2n+1}\mid i\leq 2n+1}.\qedhere\]
\end{proof}
As a consequence, we obtain a generalization of a result of \cite{ContiSalamon} concerning the intrinsic torsion of an $\SU(n)$-structure:
\begin{proposition}
\label{prop:SUnIsStronglyAdmissible}
The flag $E_0\subsetneq\dotsb\subsetneq E_{2n+1}$ is $\mathcal{I}^{\SU(n)}$-ordinary, and the group $\SU(n)\subset \LieG{O}(2n+1)$ is strongly admissible. In particular, given an $\SU(n)$-structure $(\alpha,F,\Omega)$ on $M^{2n+1}$, its intrinsic torsion is entirely determined by $d\alpha,dF,d\Omega$.
\end{proposition}
\begin{proof}
By Cartan's test and Proposition~\ref{prop:SpecialGeometries},
\[\dim T\otimes\su(n)^\perp =2n^3+3n^2+3n+1\leq \codim Z^{\SU(n)}\leq \dim T\otimes\su(n)^\perp,\]
so equality must hold and the statement follows.
\end{proof}

\section{Hypo and nearly hypo evolution equations}
\label{sec:hypo}
In this section we come to a concrete application of Theorem~\ref{thm:AbstractEmbedding} concerning the group $\SU(n)$. Let $T=\R^{2n}$,  $n\geq 2$, with basis $e_1,\dotsc,e_{2n}$, and let $F$, $\Omega$ be as in \eqref{eqn:SUnForms}.

The group $\SU(n)$ is both admissible, meaning that $\SU(n)$ is the subgroup of $\GL(n,\R)$ that fixes $(\Lambda^*T)^{\SU(n)}$, and strongly admissible, meaning that integrals of $\mathcal{I}_f$ are structures with constant intrinsic torsion. However, the constant intrinsic torsion geometries  that occur in this case are essentially only two.
\begin{proposition}
\label{prop:derivation}
 Let $f$ be a differential operator on  $(\Lambda^*T)^{\SU(n)}$ that extends a derivation of degree one on $\Lambda^*T$. Then either $f=0$ or $n=3$ and
\[f(F)=\lambda \Omega^+ + \mu \Omega^-,\quad f(\Omega^+)=\frac23\mu F^2,\quad f(\Omega^-)=-\frac23\lambda F^2.\]
\end{proposition}
\begin{proof}
If $n\neq 3$ there are no invariant three-forms; thus, $f(F)=0$. Suppose $f\neq 0$.  Then the space of invariant $n+1$-forms must be nontrivial, which implies that $n$ is odd and $f(\Omega^\pm)$ is a multiple of $F^{\frac{n+1}2}$. By the Leibnitz rule, it follows that
\[0=f(\Omega^\pm\wedge F)=f(\Omega^\pm)\wedge F,\]
and so $f(\Omega^\pm)=0$.

For $n=3$, it follows from the Leibnitz rule and $F^3=\frac32\Omega^+\wedge\Omega^-$ that $f$ must be like in the statement. Any such $f$ extends to a derivation of $\Lambda^*T$ by setting
\[f(\alpha)=\alpha\hook(\lambda\Omega^+ + \mu\Omega^-), \quad \alpha\in T.\qedhere\]
\end{proof}
So, for the group $\SU(n)$ acting on $\R^{2n}$ there are only two geometries to consider in our setup, namely Calabi-Yau geometry ($f=0$) and nearly-K\"ahler geometry ($n=3$, $f\neq 0$). In the latter case, it is customary to normalize the constants so that $\lambda=3$, $\mu=0$, although the normalization is irrelevant to the discussion to follow.

Now take  $W\subset T$ of codimension one. Since $\SU(n)$ acts transitively on the sphere in $T$, all choices of $W$ are equivalent; we shall fix
\[W=\Span{e^1,\dotsc, e^{2n-1}}.\]
Then the induced structure group $H$ is
\[\SU(n-1)=\SU(n)\cap \LieG{O}({2n-1}).\]
The space $(\Lambda^*W)^{\SU(n)}$ is spanned by the forms $\alpha$, $F$ and $\Omega^\pm$ defiend in \eqref{eqn:oddSUnForms}, with a shift in the value of $n$. On the other hand $p(\Lambda^*T(\SU(n)))$ is spanned by the forms
\[\alpha\wedge\Omega^+, \quad \alpha\wedge\Omega^-, F.\]
Thus the induced differential operator $0_W$ satisfies
\[0_W(F)=0,\quad 0_W(\alpha\wedge\Omega^\pm)=0\;.\]
Since all codimension one subspaces of $T$ are conjugate under $\SU(n)$, given a manifold $M$ with a $\SU(n)$-structure $P_{\SU(n)}$, any hypersurface $N\subset M$ admits a $\SU(n-1)$\nobreakdash-structure $P_{\SU(n-1)}$ such that $(N,P_{\SU(n-1)})$ is embedded in $(M,P_{\SU(n)})$ with type $W\subset T$. Thus, Proposition~\ref{prop:submanifold} reduces to the known fact that an oriented hypersurface in a manifold with holonomy $\SU(n)$ admits a $\SU(n-1)$-structure which is an integral of $\mathcal{I}_{0_W}$, called a \dfn{hypo} structure (see \cite{ContiSalamon, ContiFino}). On the other hand if $f$ is the nearly-K\"ahler differential operator, the induced differential operator $f_W$ is characterized by
\[f_W(F)=2\alpha\wedge\Omega^+,\quad f_W(\alpha\wedge\Omega^-)=3F\wedge F\;.\]
Again by Proposition~\ref{prop:submanifold}, we recover the known fact that  an oriented hypersurface in a nearly-K\"ahler $6$-manifold admits an $\SU(2)$-structure which is an integral of $\mathcal{I}_{f_W}$, called a \dfn{nearly-hypo} structure (see \cite{Fernandez:NearlyHypo}).

In either case, if the hypersurface $N\subset M$ is compact, one can use the exponential map to identify a tubular neighbourhood of $N$ with $N\times(a,b)$, and rewrite the  K\"ahler form and complex volume on $M$ as
\[\alpha(t)\wedge dt + F(t),\quad (\alpha(t)+idt)\wedge\Omega(t)\]
where $t$ is a coordinate on $(a,b)$ and $(\alpha(t),\Omega(t),F(t))$ is a one-parameter family of $\SU(n)$-structures on $N$. In this language, finding an embedding of $(N,P_{\SU(n-1)})$ in $(M,P_{\SU(n)})$  amounts to finding a solution of certain evolution equations (see \cite{ContiSalamon,ContiFino,Fernandez:NearlyHypo}).

A sketch of the proof of the existence of a solution to the evolution equations was given in \cite{ContiFino}, while \cite{Fernandez:NearlyHypo} left it as an open problem.  We can now give a complete and simultaneous proof of both cases, in the guise of the following:
\begin{theorem}
\label{thm:nearlyhypo}
If $(N,\alpha,F,\Omega)$ is a real analytic hypo compact manifold of dimension $2n-1$, the  hypo evolution equations admit a solution, and determine an embedding $i\colon N\to M$ into a real analytic $2n$-manifold with holonomy $\SU(n)$. If $(N,\alpha,F,\Omega)$ is a real analytic nearly-hypo compact $5$-manifold, the nearly-hypo evolution equations admit a solution, and determine an embedding $i\colon N\to M$ into a real analytic manifold with a nearly-K\"ahler $\SU(3)$-structure.
\end{theorem}
\begin{proof}
Since $\SU(n)$ is strongly admissible,  Proposition~\ref{prop:cit} implies that
\[\codim V_n(\mathcal{I}_f)=\dim \R^{2n}\otimes\frac{\so(2n)}{\su(n)}.\]
Then the flag $E_0\subsetneq\dotsb\subsetneq E_{2n}=T$ is $\mathcal{I}^{\SU(n)}$-ordinary, implying that $\R^{2n-1}\subset \R^{2n}$ is relatively admissible. On the other hand the differential operator $f$  extends to a derivation of $\Lambda^*T$ by Proposition~\ref{prop:derivation}. The statement now follows from Theorem~\ref{thm:AbstractEmbedding}.
\end{proof}
By the same argument, we also obtain a complete proof of a result stated in \cite{Bryant:Calibrated}:
\begin{corollary}
Every real analytic, parallelizable, compact Riemannian $n$-manifold can be embedded isometrically as a special Lagrangian submanifold in a manifold with holonomy $\SU(n)$.
\end{corollary}

Notice that the assumption of real analyticity refers not only to the manifold, but to the structure as well.

\section{$\alpha$-Einstein-Sasaki geometry and hypersurfaces}
\label{sec:go}
In this section we classify the constant intrinsic torsion geometries for the group $\SU(n)\subset \LieG{O}(2n+1)$, and write down evolution equations for hypersurfaces which are orthogonal to the characteristic direction, in analogy with  Section~\ref{sec:hypo}.

Let $T=\R^{2n+1}$. The space $(\Lambda^*T)^{\SU(n)}$ is spanned by the forms $\alpha$, $F$ and $\Omega^\pm$, defined in \eqref{eqn:oddSUnForms}. In order to classify the differential operators on $(\Lambda^*T)^{\SU(n)}$, we observe that every element $g$ of the normalizer $N(\SU(n))$ of $\SU(n)$ in \mbox{$\SO(2n+1)$} maps $(\Lambda^*T)^{\SU(n)}$ to itself; this defines a natural notion of equivalence among differential operators.
\begin{proposition}
\label{prop:oddderivation}
Let $f$ be a derivation of $(\Lambda^*T)^{\SU(n)}$; then $f$ is a differential operator that extend to a derivation of degree one on $\Lambda^*T$ if and only if one of the following holds:
\begin{itemize}
\item[(A)] $f(\alpha)=0$, $f(F)=2\lambda\alpha\wedge F$, $f(\Omega)=n(\lambda-\mu i)\alpha\wedge\Omega$;
\item[(B)] $f(\alpha)=\lambda F$, $f(F)=0$, $f(\Omega)=-\mu i\alpha\wedge\Omega$;
\item[(C)] $n=2$, and up to $N(\SU(2))$ action, $\tilde f$ has the form (A) or (B);
\item[(D)] $n=3$, $f(\alpha)=0$, $f(F)=3\lambda\Omega^- - 3\mu\Omega^+$, $f(\Omega)=2(\lambda+i\mu)F^2$;
\end{itemize}
here $\lambda$ and $\mu$ are real constants.
\end{proposition}
\begin{proof}
By Lemma~\ref{lemma:ZEmpty}, we have to classify the $\SU(n)$-equivariant derivations of degree one $\tilde f$ of $\Lambda^*T$ with $\tilde f^2=0$ on $(\Lambda^*T)^{\SU(n)}$. As a representation of $\SU(n)$, we have
\[\Hom(T,\Lambda^2T)=(\R\oplus[\Lambda^{1,0}])\otimes ([\Lambda^{2,0}]\oplus \leftdoublebracket\Lambda^{1,1}_0\rightdoublebracket\oplus\R\oplus [\Lambda^{1,0}]).\]
Decomposing into irreducible components, this tensor product always contains three trivial components, one of which is in $\R\otimes\R$, and the other two in $[\Lambda^{1,0}]\otimes [\Lambda^{1,0}]$. Moreover if $n=2$ both $[\Lambda^{1,0}]\otimes [\Lambda^{1,0}]$ and $\R\otimes[\Lambda^{2,0}]$ contain each two more trivial components, and if $n=3$ $[\Lambda^{1,0}]\otimes [\Lambda^{2,0}]$ contains two trivial components. Thus, we have to consider three different cases.

(\emph{i}) If $n>3$ and $\tilde f$ is an invariant derivation, there are constants $\lambda$, $\mu$ and $k$ such that
\[\tilde f(\alpha) = kF, \quad \tilde f(e^i)=\alpha\wedge (\lambda e^i+\mu e_i\hook F),\]
so
$\tilde f(F)=2\lambda \alpha \wedge F$. Since $f$ is a differential operator, it follows that
\[0=\tilde f(\lambda\alpha\wedge F)=\lambda k F\wedge F,\]
hence either $\lambda=0$ or $k=0$. On the other hand, in complex terms,
\[\tilde{f}(\Omega)=n(\lambda-\mu i)\alpha\wedge\Omega\;,\]
so we obtain (A) or (B). Notice that these are indeed differential operators, i.e. $f^2=0$.

(\emph{ii}) If $n=3$, the invariant derivation $\tilde f$ has the form
\[\tilde f(\alpha) = kF, \quad \tilde f(e^i)=\alpha\wedge (\lambda e^i+\mu e_i\hook F)+\beta e_i\hook\Omega^+ + \gamma e_i\hook\Omega^-,\]
and therefore
\[\tilde f(F)=2\lambda \alpha \wedge F +3\beta\Omega^- - 3\gamma\Omega^+,\quad
\tilde f(\Omega)=2(\beta+i\gamma)F^2+3(\lambda-i\mu)\alpha\wedge\Omega.\]
If $k\neq 0$, then
\[0=\tilde f^2(\alpha)=2\lambda \alpha \wedge F +3\beta\Omega^- - 3\gamma\Omega^+,\]
so $\lambda=0=\beta=\gamma$ and $\tilde f$ lies in the second family. On the other hand if $k=0$, then
 \[\tilde f(2\lambda \alpha \wedge F +3\beta\Omega^- - 3\gamma\Omega^+)=3(\lambda\beta-3\mu\gamma)\alpha\wedge\Omega^- -3 (\lambda\gamma+3\mu\beta)\alpha\wedge\Omega^+\]
must be zero. Similarly,
\[0=\tilde f^2(\Omega)=(\lambda+3i\mu)(\beta+i\gamma)\alpha \wedge F^2,\]
so either $\beta=\gamma=0$, or $\lambda=\mu=0$, corresponding to (A) and (D) in the statement.

(\emph{iii}) If $n=2$,  $\tilde f$ has the form
\[\tilde f(\alpha) = kF+\beta\Omega^+ + \gamma\Omega^-, \quad \tilde f(e^i)=\alpha\wedge (\lambda e^i+\mu e_i\hook F +\sigma e_i\hook\Omega^++\tau e_i\hook\Omega^-)\;;\]
since the normalizer $N(\SU(2))$ contains a copy of $\SU(2)$ that acts as rotations on the space spanned by $F$ and $\Omega^\pm$, we can assume that
$\mu  F +\sigma \Omega^++\tau \Omega^-$ is a multiple of $F$, i.e. $\sigma=0=\tau$. We obtain
\[
\tilde f(F)=2\lambda F\wedge\alpha,\quad
\tilde f(\Omega)=2(\lambda-\mu i)\alpha\wedge\Omega
\]
The condition  $f^2=0$ is then equivalent to
\begin{gather*}
 \lambda k=0, \quad \lambda\beta=0, \quad  \lambda\gamma =0, \quad \mu\gamma=0, \quad \mu\beta=0.
\end{gather*}
Now if $\gamma=0=\beta$ we are in case (A) or (B). Otherwise, we have $\mu=0$, i.e.
\[\tilde f(\omega)=2\lambda \omega\wedge\alpha, \quad \omega\in\Span{F,\Omega^+,\Omega^-}.\]
Thus, we can use the action of $N(\SU(2))$ to reduce to the case where $\tilde{f}(\alpha)$ is a multiple of $F$.
\end{proof}
The differential operators appearing in Proposition~\ref{prop:oddderivation}, to be denoted each by the corresponding letter, can be intepreted as follows.
The differential operator $A$ determines a codimension one foliation $\ker\alpha$, where each leaf has an integrable induced $\SU(n)$-structure. Conversely, given an integrable $\SU(n)$-structure $(F,\Omega)$ on $M^{2n}$, we can define an $\SU(n)$-structure on $M^{2n}\times\R$ by
\[\tilde F = e^{2\lambda t} F,\quad  \alpha = dt, \quad \tilde\Omega = e^{n(\lambda - i\mu) t}\Omega,\]
which is an integral of $\mathcal{I}_A$.

Similarly, the differential operator $C$ corresponds to a foliation with nearly-K\"ahler leaves. On the other hand $B$ corresponds to $\alpha$-Einstein-Sasaki geometry  (see \cite{ContiSalamon} for a proof of this fact in the five-dimensional case). The methods of Section~\ref{sec:hypo} apply with minimal changes; the main difference is that here $\SU(n)$ does not act transitively on the sphere in $T=\R^{2n+1}$. Indeed, codimension one subspaces $W\subset T$ have an invariant, namely the angle $\Gamma$ that they form with the characteristic direction $e_{2n+1}$. However, by \cite{Conti:CohomogeneityOne}, if $T$ is $\alpha$-Einstein-Sasaki and one uses the exponential to identify a tubular neighbourhood of a hypersurface $N\subset T$ with the product $M\times(a,b)$, the angle $\Gamma$ is constant along the radial direction.

Let us consider the case that $W$ is tangent to the characteristic direction, i.e.
\[W=\Span{e_1,\dots,e_{2n-1},e_{2n+1}}.\]
Then $\SU(n)\cap \LieG{O}(W)=\SU(n-1)$, and  the space $(\Lambda^*W)^{\SU(n-1)}$ is spanned by
\[\alpha=e^{2n+1},\quad \beta = e^{2n-1}, \quad F=e^{12}+\dotsb+e^{2n-3,2n-2},\]
and the real and imaginary part of
\[\Omega^+ + i\Omega^-=(e^1+ie^2)\dotsm (e^{2n-3}+ie^{2n-2}).\]
Then  $p(\Lambda^*T(\SU(n)))$ is spanned by the forms $\alpha$, $F$ and $\Omega\wedge\beta$; the projection induces an operator satisfying
\begin{equation}
 \label{eqn:Go}
B_W(\alpha)=\lambda F, \quad B_W(\Omega\wedge\beta)=-i\mu\alpha\wedge\Omega\wedge\beta, \quad B_W(F)=0.
\end{equation}

Given an integral $(\alpha(0),F(0),\Omega(0),\beta(0))$ of $\mathcal{I}_{B_W}$ on a $2n$-dimensional manifold $N$, embedding $N$ into a manifold $M$ with an $\SU(n)$-structure which is an integral of $\mathcal{I}_B$  is equivalent to extending $(\alpha(0),F(0),\Omega(0),\beta(0))$ to a one-parameter family of $\SU(n-1)$-structures $(\alpha(t),F(t),\Omega(t),\beta(t))$ such that
\begin{equation}
\label{eqn:GoEvolution}
 \frac{\partial }{\partial t}\alpha = -\lambda \beta, \quad  \frac{\partial }{\partial t}F = -d\beta, \quad
 \frac{\partial }{\partial t}(\beta\wedge\Omega) = id\Omega -\mu\alpha\wedge\Omega
\end{equation}
Indeed, in that case $(\alpha(t),F(t),\Omega(t),\beta(t))$ is an integral of $\mathcal{I}_{B_W}$ for all $t$, and the forms
\[\alpha(t), \quad \Omega(t)\wedge (\beta(t)+idt), \quad F(t)+\beta(t)\wedge dt\]
define a $\SU(n)$-structure on $N\times(a,b)$ which is an integral of $\mathcal{I}_B$. The converse follows from the fact that the angle $\Gamma$ is constant along the radial direction.

Applying Theorem~\ref{thm:AbstractEmbedding}, we obtain an odd-dimensional (with reference to the ambient space) analogue of Theorem~\ref{thm:nearlyhypo}.
\begin{theorem}
\label{thm:go}
If $N$ is a real analytic compact manifold of dimension $2n$ with a real analytic $\SU(n-1)$-structure $(\alpha,F,\Omega,\beta)$ which is an integral of $\mathcal{I}_{B_W}$ as defined in \eqref{eqn:Go}, the evolution equations \eqref{eqn:GoEvolution} admit a solution with initial data $(\alpha,F,\Omega,\beta)$, and determine an embedding $i\colon N\to M$ into a real analytic $\alpha$-Einstein-Sasaki manifold as a hypersurface tangent to the characteristic direction.
\end{theorem}

\section{An irreducible geometry in dimension $9$}
\label{sec:so3}
As a final application, we consider a special geometry modelled on the $9$\nobreakdash-dimensional irreducible representation $T$ of $\SO(3)$. This geometry is described by a $4$-form and a $5$-form; if closed, the $5$-form defines a calibration, and one can look for calibrated embeddings of five-manifolds in this sense; more generally, one can fix a five-dimensional subspace of $\R^9$, and look for ``compatible'' embeddings.

To make things explicit, we fix a basis of $\so(3)$ satisfying
\[[H,X]=\sqrt2Y, \quad [H,Y]=-\sqrt2X, \quad [X,Y]=\sqrt2H\;;\]
and identify the representation $T$ with $\R^9$ in such a way that the generic element $aH+bX+cY$ of $\so(3)$ acts as the matrix
\[  \left(\begin{array}{ccccccccc}0&2  b&0&0&-4  a \sqrt{2}&-2  c&0&0&0\\-2  b&0& \sqrt{7} b&0&-2  c&-3  a \sqrt{2}&- c \sqrt{7}&0&0\\0&- \sqrt{7} b&0&3  b&0&- c \sqrt{7}&-2  a \sqrt{2}&-3  c&0\\0&0&-3  b&0&0&0&-3  c&- a \sqrt{2}&2  \sqrt{5} b\\4  a \sqrt{2}&2  c&0&0&0&2  b&0&0&0\\2  c&3  a \sqrt{2}& c \sqrt{7}&0&-2  b&0& \sqrt{7} b&0&0\\0& c \sqrt{7}&2  a \sqrt{2}&3  c&0&- \sqrt{7} b&0&3  b&0\\0&0&3  c& a \sqrt{2}&0&0&-3  b&0&2  c \sqrt{5}\\0&0&0&-2  \sqrt{5} b&0&0&0&-2  c \sqrt{5}&0\end{array}\right)
\]
Since this matrix is skew-symmetric, the  action of $\SO(3)$ preserves the standard metric on $\R^9$.

By standard representation theory, one sees that the action of $\SO(3)$ on  $\Lambda^*T$ leaves invariant a $4$-form $\gamma$, and this is the only invariant along with its Hodge dual $*\gamma$ and the volume form. If $e^1,\dotsc, e^9$ is the standard basis of $\R^9$, using the expression for the $\SU(2)$ action given above we can identify $\gamma$ as
\begin{multline}
\label{eqn:SO3gamma}
\gamma=\frac{1}{4}\sqrt5\left( -e^{2589}+ e^{1249} -e^{1689}+ e^{4569} \right)
- (e^{1357}+e^{1256})
+\frac{7}{8}  e^{3478}\\
+\frac{1}{8}\sqrt{35}\left(-e^{3689} -e^{2789}+e^{4679}+ e^{2349} \right)
-\frac{1}{2} e^{1458}
+\frac{1}{8}  e^{2367}\\
+\frac38(e^{3456}+e^{5678}-e^{2468}-e^{2457}-e^{2358}-e^{1467}-e^{1368}+e^{1278}+e^{1234})
\end{multline}
Thus, on a $9$-manifold $M$ with an $\SO(3)$-structure one has two canonical forms, which we also denote by  $\gamma$ and $*\gamma$. The intrinsic torsion space has dimension $279$ and splits into $25$ irreducible components, but the most natural torsion classes to consider are those defined by one of
\begin{equation}
 \label{eqn:SO3Geometries} d\gamma=0, \quad d*\gamma=0, \quad d\gamma=0=d*\gamma.
\end{equation}
 By Section~\ref{sec:Preparatory}, to each geometry in \eqref{eqn:SO3Geometries} one can associate an exterior differential system, which is involutive in the second case (by Lemma~\ref{lemma:gammaEstable} and Proposition~\ref{prop:stable}), though it does not appear to be involutive in the other two cases.

\begin{remark}
It would also be natural to consider the geometry defined by $d\gamma=\lambda *\gamma$ with $\lambda$ a non-zero constant, but structures of this type do not exist. Indeed,  let \mbox{$A=(\Lambda^*T)^{\SO(3)}$}. Then $A$ is generated by two elements $\gamma$, $*\gamma$, but the differential operator defined by $f(\gamma)=*\gamma$ does not extend to a derivation of $\Lambda^*T$, since $\Hom(T,\Lambda^2T)$ has no trivial submodule (see Lemma~\ref{lemma:ZEmpty}).
\end{remark}
Now suppose that $*\gamma$ is closed; this means that $8$ out of the $25$ components of the intrinsic torsion vanish (counting dimensions, $84$ out of $279$). Like in Section~\ref{sec:psu3}, we can renormalize $*\gamma$ into a calibration and try to obtain an embedding result. However, this time the form is not stable, but a weaker result holds:
\begin{lemma}
\label{lemma:gammaEstable}
Let $\alpha$ be a vector in the finite set
\begin{equation}
 \label{eqn:finiteset}\{e^1,e^2,e^4,e^5,e^6,e^8\}\subset\R^9.
\end{equation}
Then the form $*\gamma$ is $\alpha^\perp$-stable.
\end{lemma}
In order to prove an analogue of Theorem~\ref{thm:psu3Embedding2}, one would have to show that $\R^9$ contains calibrated subspaces orthogonal to one of \eqref{eqn:finiteset}. This appears to be a difficult problem, because the structure group is small and $*\gamma$ has a complicated form.

At any rate, one can fix a $5$-dimensional $W\subset\R^9$ and consider submanifolds $N\subset M$ embedded with type $W\subset\R^9$, in the sense that at each point $x\in N$ one can choose a coframe $e^1,\dotsc, e^9$ such that $\gamma$ has the form \eqref{eqn:SO3gamma} and the coframe maps $T_xN\subset T_xM$ to $W\subset\R^9$.
\begin{theorem}
\label{thm:so3Embedding}
If $N$ is a compact, parallelizable, real analytic Riemannian $5$-manifold and $W\subset\R^9$ is a $5$-dimensional subspace orthogonal to one of the vectors \[e^1,e^2,e^4,e^5,e^6,e^8\]
then $N$ can be embedded with type $W\subset\R^9$ in a $9$-manifold $M$ with an $\SO(3)$-structure such that $*\gamma$ is closed.
\end{theorem}
\begin{proof}
 Lemma~\ref{lemma:gammaEstable} identifies six $8$\nobreakdash-dimensional subspaces $E\subset\R^9$ such that the $5$\nobreakdash-form $*\gamma$ is $E$-stable. By Proposition~\ref{prop:stable}, if $W$ is contained in such an $E$, then $W\subset \R^9$ is relatively admissible, and by Theorem~\ref{thm:AbstractEmbedding} the statement follows.
\end{proof}
We leave it as an open problem to determine the subspaces $W\subset\R^9$ calibrated by $*\gamma$.

\bibliographystyle{plain}
\bibliography{torsion}

\end{document}